\documentclass[twocolumn,10pt]{autart}    

\usepackage{graphicx}          
\usepackage{algorithm}            
\usepackage{amsmath}
\usepackage{enumerate}
\usepackage{natbib}                          

\usepackage{amsfonts}
\usepackage{algpseudocode}
\usepackage{amssymb}
\usepackage{subfigure}
\usepackage{hyperref}
\hypersetup{
    colorlinks = true,
    citecolor = cyan,
}
\usepackage{wrapfig}

\newtheorem{remark}{Remark}[section]

\newtheorem{assumption}{Assumption}[section]

\newtheorem{lemmax}{Lemma}[section]

\newtheorem{theorem}{Theorem}[section]

\begin{document}

\begin{frontmatter}

\title{A Quantized Order Estimator} 

\author{Lida~Jing}\ead{jing@sdu.edu.cn}


\address{School of Mathematics, Shandong University, Jinan, Shandong 250100, China}


\begin{abstract}                          
This paper considers the order estimation problem of stochastic autoregressive exogenous input (ARX) systems by using quantized data. Based on the least squares algorithm and inspired by the control systems information criterion (CIC), \textcolor{blue}{a new kind of criterion aimed at addressing the inaccuracy of quantized data is} proposed for ARX systems with quantized data. When the upper bounds of the system orders are known and the persistent excitation condition is satisfied, the system order estimates are shown to be consistent for small quantization step. Furthermore, a concrete method is given for choosing quantization parameters to ensure that the system order estimates are consistent. A numerical example is given to demonstrate the effectiveness of the theoretical results of the paper.
\end{abstract}

\end{frontmatter}

\section{Introduction}
System identification with quantized data is a challenging research topic \citep{Wang2003,Gustafsson2009}. \textcolor{blue}{In many cases, using quantized data during the system identification process will bring quantization error, which increases the difficulty of analysis.} Up to now, a large number of identification methods with quantized data have been developed, including \citep{Wang2003,Wang2010,Jing2019,Jing2021,Jing2022,Wang2019,Zhang2019,Diao2020}, to name a few. In particular, \citep{Wang2003} proposed two different frameworks, namely, stochastic and deterministic frameworks so as to identify systems. \citep{Wang2010} gave some motivating examples of quantized measurements and introduced the methods and algorithms of system identification for set-valued linear systems. \citep{Jing2019} used projection algorithm to estimate parameters of quantized deterministic autoregressive moving average (DARMA) systems, and proved the boundedness of parameter estimation error by designing system inputs. \citep{Wang2019} researched the identification of multi-agent systems with quantized observations. \citep{Zhang2019} concerned the system identification for FIR systems with set-valued \textcolor{blue}{and} precise data received from multiple sensors. \citep{Jing2021,Jing2022} solved the parameter estimation problem of quantized DARMA systems and quantized stochastic autoregressive exogenous input (ARX) systems with the help of the least squares, respectively.

The system identification task for ARX systems consists of estimating (i) the orders, (ii) the parameters, and (iii) the covariance matrix of system noise. However, the contributions listed above are all for parameter estimation with quantized data. As for order estimation by using quantized data, it is a novel problem. Obviously, selecting the right model order is the first step for the goal of estimating system parameters. A number of classic order estimation techniques such as \citep{Akaike1969,Soderstrom1977,Hannan1980,Soderstrom1989,Liang1993,Hannan1979,Hannan1982,Chen1987,Guo1989} have been made since about the 1970s. Specifically, Akaike proposed a well-known criterion, Akaike's Information Criterion (AIC) \citep{Akaike1969}. \citep{Soderstrom1977} proved that Final Prediction-Error (FPE) criterion and AIC are asymptotically equivalent. \citep{Hannan1979} proved that a strongly consistent estimation \textcolor{blue}{of the order can be based on} the law of iterated logarithm for the partial autocorrelations. \citep{Hannan1980} made some consistent works on the order estimation. \citep{Hannan1982} established the asymptotic properties under very general conditions. \citep{Chen1987} got a consistent estimate of the order of feedback control systems with system parameters estimated by the least squares method. \citep{Guo1989} introduced a new criterion, control systems information criterion (CIC), so as to estimate orders of the linear stochastic feedback control system. \citep{Liang1993} proposed an approach for model order determination based on the minimum description length (MDL) criterion which is shown to depend on the minimum eigenvalues of a covariance matrix derived from the observed data. 

Considering the wide use of quantized data and the important value of order estimation, it is of significance to study order estimation based on quantized data. The introduction of quantized data will produce quantization error, which brings difficulties to order estimation. By using some conclusions of \citep{Jing2022}, one order estimation method of ARX models with uniform quantized data is proposed. \textcolor{blue}{The order estimation algorithm in the paper is utilized in the following process. First of all, the range of ARX system orders is selected (i.e., $0\leq p\leq p_{max}$ and $0\leq q\leq q_{max}$, where $p$ is the order of the AR part and $q$ is the order of the exogenous part). Then for each $(p, q)$ pair the parameters of the model are
estimated by the least squares} under the assumption that $p$ and $q$ are the right model orders. Finally, a prediction error variance for the model is calculated by \textcolor{blue}{the proposed criterion} and the $(p, q)$ pair yielding the lowest value is chosen as the best estimate of the model order. So, the key step of estimation lies in two aspects: \textcolor{blue}{the design of a criterion for the order estimate algorithm as well as the choice of a quantization step}. In fact, they are complementary.

In contrast to the previous works \citep{Hannan1980,Liang1993,Hannan1979,Hannan1982,Chen1987,Jing2022,Wang2019}, the main contributions of this paper are summarized as follows.
\begin{itemize}
\item As mentioned earlier, order estimation is one component of system identification problems. However, to the best of my knowledge, the \textcolor{blue}{existing papers} of quantized system identification mainly focus on quantized parameter estimation. The \textcolor{blue}{discussion} about quantized order estimation is pretty rare. Actually, literatures like \citep{Jing2019,Jing2021,Jing2022,Wang2019,Zhang2019} considered quantized parameter estimation based on known system orders. And different from them, in this paper, we study the quantized order estimation problem when the system orders and parameters are both unknown.

\item Compared with classic papers \citep{Hannan1980,Soderstrom1989,Liang1993,Hannan1979,Hannan1982,Chen1987,Guo1989} on order estimation based on accurate data, we study order estimate problem under uniform quantized observations. To be more concrete, \textcolor{blue}{one} of the difficulties in \textcolor{blue}{designing order estimate algorithm} is how to make full use of the roughness of quantized observations. Quantized data make the structure of classic estimation algorithms more complex and the estimated parameter can not converge to real value in many cases. By designing \textcolor{blue}{the criterion} and using some hypotheses \textcolor{blue}{of system parameters and orders}, the quantized order estimation can converge to real value in some sense.

\item Different from \citep{Jing2019,Jing2021,Wang2019}, the model researched in this paper contains stochastic noises. So, the algorithm analysis methods in parameter estimation part of this note are quite different.
\end{itemize}

In this paper, $\mathbb{R}$ denotes real number field. For a given vector or matrix $x$, $x^{\top}$ denotes the transpose of $x$; $||x||$ denotes the Euclidean norm for vector case and the corresponding induced norm for matrix case. $\lambda_{min}\left(\right)$ denotes the smallest eigenvalue of the matrix between round brackets. The rest of the paper is as follows. In section 2, we describe the model. Section 3 shows the specific order estimation algorithm for the quantized ARX model, and the influence of quantization error on the order estimation is analyzed. Section 4 uses a numerical example to demonstrate the main result. \textcolor{blue}{Section 5 concludes this work}.

\section{Model}
Consider the following ARX \textcolor{blue}{system}:
\begin{equation}\label{(1)}
A(z)y_{n+1}=B(z)u_{n}+w_{n+1}, \quad n\geq0,
\end{equation}
where $y_{n}$, $u_{n}$ and $w_{n}$ are the system output, system input and system noise. Besides, \textcolor{blue}{define that $N(0,1)$ indicates a Gaussian distribution
with zero mean and variance $1$}. The noise \{$w_n$\} is a sequence of independent and identically distributed (i.i.d.) \textcolor{blue}{random} variables and $w_n\sim N(0,1)$. For simplicity, suppose $y_n=u_n=w_n=0$, $\forall n<0$.
\begin{equation*}
A(z)=1+a_1z+a_2z^2+\cdots+a_{p_0}z^{p_0}, \quad p_0\geq0,
\end{equation*}
\begin{equation*}
B(z)=b_1+b_2z+\cdots+b_{q_0}z^{q_0-1},\quad  q_0\geq1,
\end{equation*}
where $a_i$ and $b_j$ are unknown system parameters. $z$ is the shift-back operator and the orders $p_0$, $q_0$ are unknown. $a_{p_0}\neq0$, $b_{q_0}\neq0$.

For the convenience of proving, the model (\ref{(1)}) can be rewritten as follows:
\begin{equation}\label{(2)}
y_{n+1}=\theta^{\top}(p_0,q_0)\varphi_n(p_0,q_0)+w_{n+1},
\end{equation}
where $\theta(p_0,q_0)=\left[-a_1, \cdots, -a_{p_0}, b_1, \cdots, b_{q_0}\right]^{\top}$ and $\varphi_n(p_0,q_0)=\left[y_{n}, \cdots, y_{n-p_0+1}, u_{n}, \cdots, u_{n-q_0+1}\right]^{\top}$.

This paper considers the condition that the system output $y_{n}$ cannot be directly measured and only its quantized value is known. We want to design an order estimation algorithm and analyze the influence of the quantization step on order estimation.

For a given constant $\varepsilon>0$ and any $n$=1, 2, ... , the quantized value of $y_{n}$ is from the following uniform quantizer:
\begin{align}\label{(3)}
s_n=\varepsilon\left\lfloor\frac{y_n}{\varepsilon}+\frac{1}{2}\right\rfloor.
\end{align}
We can call $\varepsilon$ the quantization step and $s_n$ is the quantized output.

\begin{remark}
The more direct form of the equation (\ref{(3)}) is
\begin{align*}
s_{n}=\left\{
\begin{aligned}
&\vdots\\
&-2\varepsilon, &&\text{}y_{n}\in\left[-\frac{5\varepsilon}{2}, -\frac{3\varepsilon}{2}\right),\\
&-\varepsilon, &&\text{}y_{n}\in\left[-\frac{3\varepsilon}{2}, -\frac{\varepsilon}{2}\right),\\
&0,           &&\text{}y_{n}\in\left[-\frac{\varepsilon}{2}, \frac{\varepsilon}{2}\right),\\
&\varepsilon, &&\text{}y_{n}\in\left[\frac{\varepsilon}{2}, \frac{3\varepsilon}{2}\right),\\
&2\varepsilon, &&\text{}y_{n}\in\left[\frac{3\varepsilon}{2}, \frac{5\varepsilon}{2}\right),\\
&\vdots
\end{aligned}\right..
\end{align*}
\end{remark}

From (\ref{(2)}) and (\ref{(3)}) we know that
\begin{align}\label{(4)}
s_{n+1}=\theta^{\top}(p_0,q_0)\psi_n(p_0, q_0)+w_{n+1}+\epsilon_{n+1},
\end{align}
where
\begin{align}\label{(5)}
\psi_n(p_0, q_0)=\left[s_{n}, \cdots, s_{n-p_0+1}, u_{n}, \cdots, u_{n-q_0+1}\right]^{\top},
\end{align}
and $\epsilon_{n+1}$ is the quantization noise at time $n+1$, which is produced by quantized outputs and its concrete property is as follows.

From (\ref{(2)}), (\ref{(4)}) we know that
\begin{align}\label{(6)}
\left|\epsilon_{n+1}\right|=&\left|s_{n+1}-\theta^{\top}(p_0,q_0)\psi_n(p_0, q_0)-w_{n+1}\right|\notag\\
=&\left|s_{n+1}-\theta^{\top}(p_0,q_0)\psi_n(p_0, q_0)\right.\notag\\
&\left.-\left(y_{n+1}-\theta^{\top}(p_0,q_0)\varphi_n(p_0,q_0)\right)\right|\notag\\
=&\left|s_{n+1}-y_{n+1}\right.\notag\\
&\left.+\theta^{\top}(p_0,q_0)\left(\varphi_n(p_0,q_0)-\psi_n(p_0, q_0)\right)\right|\notag\\
\leq&\left|s_{n+1}-y_{n+1}\right|\notag\\
&+\left|\theta^{\top}(p_0,q_0)\left(\varphi_n(p_0,q_0)-\psi_n(p_0, q_0)\right)\right|\notag\\
\leq&\frac{\varepsilon}{2}+\frac{\varepsilon}{2}\left(\left|a_1\right|+\left|a_2\right|+\cdots+\left|a_{p_0}\right|\right)\notag\\
=&\frac{\varepsilon}{2}\left(\left|a_1\right|+\left|a_2\right|+\cdots+\left|a_{p_0}\right|+1\right).
\end{align}
So, we can assume $\epsilon_{n}$ is the bounded noise.

\section{Order estimation of quantized ARX systems}
The purpose of this paper is to estimate $p_0$ and $q_0$ in (\ref{(4)}) by using system inputs and quantized outputs. In this section, we give the specific order estimate method and analyze its properties.

Define
\begin{equation}\label{(7)}
\begin{cases}
\psi_i(p, q):=\left[s_{i}, \cdots, s_{i-p+1}, u_{i}, \cdots, u_{i-q+1}\right]^{\top},\\
P_{n+1}(p, q):=\left(I+\sum_{i=0}^{n}\psi_i(p, q)\psi^{\top}_i(p, q)\right)^{-1},
\end{cases}
\end{equation}
\textcolor{blue}{where $s_i=u_i=0$, when $i\leq0.$ And define} $\lambda^{(p, q)}_{min}(n)$ the smallest eigenvalue of $P^{-1}_{n+1}(p, q)$.
\subsection{Assumptions}
In order to proceed the analysis, we introduce the following assumptions.

\begin{assumption}\label{a1}
\{$u_i$\} is a sequence of independent and identically distributed (i.i.d.) \textcolor{blue}{random} variables and $u_i$ satisfies uniform distribution in $[-\delta, \delta]$, $\delta>0$.
\end{assumption}
\begin{assumption}\label{a2}
$A(z)$ is stable, i.e., $A(z)\neq0$, $\forall|z|\leq1$.
\end{assumption}
\begin{assumption}\label{a3}
There exists a constant $c>0$ such that $\left|a_i\right|\leq c$, $\left|b_j\right|\leq c$, $i=1, ..., p_0$, $j=1, ..., q_0$, and $\varepsilon<\frac{1}{2\left(1+p_0c\right)}$.
\end{assumption}
\begin{assumption}\label{a4}
$\{p_0, q_0\}$ belongs to a known finite \textcolor{blue}{set} $M$:
$$M\triangleq\left\{(p, q): 0\leq p\leq p^*, 1\leq q\leq q^*\right\},$$
where the integers $p^*>0$, $q^*>0$.
\end{assumption}

\begin{assumption}\label{a5}
There exists a constant $c_1>0$ such that
$$
\lambda^{(p, q^*)}_{min}(n)\geq c_1\left(n+1\right), a.s., n\to\infty.$$
for all $0\leq p\leq p^*$.
\end{assumption}

\begin{assumption}\label{a6}
There exists a constant $c_2>0$ such that
$$
\lambda^{(p^*, q)}_{min}(n)\geq c_2\left(n+1\right), a.s., n\to\infty.$$
for all $0\leq q\leq q^*$.
\end{assumption}

\begin{remark}
Assumption \ref{a1} means system inputs \{$u_i$\} are bounded and satisfy \textcolor{blue}{uniform distribution}. Assumptions \ref{a2} and \ref{a4} are common in classic system identification literature. Assumption \ref{a3} is always used in quantized identification. Assumptions \ref{a5} and \ref{a6} mean persistent excitation condition can be satisfied and they are pretty important to the proof of theorem in the paper.
\end{remark}

\subsection{The estimation of $p_0$}
In this section, we will prove the convergence of \textcolor{blue}{the estimate of $p_0$}.

First, we give the analyses of the matrix composed by quantized regressor vectors.
\begin{lemmax}\label{l1}
Suppose Assumptions \ref{a1}\mbox{-}\ref{a2} are satisfied. Then, as $n\to\infty$, there is a constant $c_3>0$ such that
\begin{align}\label{(11)}
\lambda^{(p_0, q^*)}_{max}(n)\leq c_3\left(n+1\right), a.s.,
\end{align}
where $\lambda^{(p_0, q^*)}_{max}(n)$ denotes the largest eigenvalue of
$
\sum_{i=0}^{n}\psi_i(p_0, q^*)\psi^{\top}_i(p_0, q^*)+I.
$
\end{lemmax}
\emph{Proof}: The proof can be seen in Appendix \ref{A}.

\textcolor{blue}{Define}
\begin{align}\label{(13)}
\bar{\theta}(p, q)=\left[-a_1, \cdots, -a_p, b_1, \cdots, b_q\right]^{\top},
\end{align}
where
\begin{align}\label{(14)}
a_i=0, b_j=0,\quad i>p_0, j>q_0.
\end{align}

\textcolor{blue}{And the estimation of $\bar{\theta}(p, q)$ is defined as}
\begin{align}\label{(15)}
\theta_n(p, q):=&\left(\sum^{n-1}_{i=0}\psi_i(p, q)\psi^{\top}_i(p, q)+I\right)^{-1}\sum^{n-1}_{i=0}\psi_i(p, q)s_{i+1}\notag\\
=&P_n(p, q)\sum^{n-1}_{i=0}\psi_i(p, q)s_{i+1},
\end{align}
where
\begin{align}\label{(16)}
\theta_n(p, q)=\left[-a_{1n}, \cdots, -a_{pn}, b_{1n}, \cdots, b_{qn}\right]^{\top}.
\end{align}

\begin{lemmax} \label{l2}
Suppose Assumptions \ref{a1}\mbox{-}\ref{a5} are satisfied. Then as $n\to\infty$,
\begin{align}\label{(17)}
&\left|\left|\left(\sum^{n-1}_{i=0}\psi_i(p_0, q^*)\psi^{\top}_i(p_0, q^*)+ I\right)^{-\frac{1}{2}}\notag\right.\right.\\
&\left.\left.\sum_{i=0}^{n-1}\psi_i(p_0, q^*)\left(w_{i+1}+\epsilon_{i+1}\right)\right|\right|^2\notag\\
\leq&\left(1+p_0c\right)\varepsilon n+o\left(n\right), a.s.
\end{align}
\end{lemmax}
\emph{Proof}: The proof can be seen in Appendix \ref{B}.

Next, we show the properties of parameter estimation error.
\begin{lemmax}\label{l3}
Suppose Assumptions \ref{a1}\mbox{-}\ref{a5} are satisfied under the condition $p\leq p_0$, and define
\begin{align}\label{(18)}
\hat{\theta}_n(p):=[&-a_{1n}(p), \cdots, -a_{pn}(p), \underbrace{0, \cdots, 0}_{p_0-p},\notag\\
&b_{1n}(p), \cdots, b_{q^*n}(p)]^{\top},
\end{align}
where $a_{in}(p)$, $b_{in}(p)$ are of $\theta_n(p, q^*)$. 

Let
\begin{align}\label{(19)}
\tilde{\theta}_n(p)=\bar{\theta}(p_0, q^*)-\hat{\theta}_n(p).
\end{align}
Then as $n\to\infty$, there is a constant $\gamma$ such that
\begin{align}\label{(20)}
\left|\left|\tilde{\theta}_n(p)\right|\right|\leq\gamma, a.s.
\end{align}
\end{lemmax}
\emph{Proof}: The proof can be seen in Appendix \ref{C}.

Then, we give the form of quantized criterion $L_n(p, q)$ and the order estimation algorithm.

Define
\begin{align}\label{(21)}
L_n(p, q):=\sigma_n(p, q)+l_n\cdot\left(p+q\right),
\end{align}
where
\begin{align}\label{(22)}
\sigma_n(p, q)=\sum^{n-1}_{i=0}\left(s_{i+1}-\theta_n^{\top}(p, q)\psi_i(p, q)\right)^2,
\end{align}
and the restrictions of $l_n$ will be given later.

The order estimation $\hat{p}_n$ of $p_0$ is defined as
\begin{align}\label{(23)}
\hat{p}_n:=\mathrm{argmin}_{\substack{0\leq p\leq p^*}}L_n(p, q^*).
\end{align}

Now, we give the upper bound of $\sigma_n(p_0, q^*)$ in the following lemma.
\begin{lemmax}\label{l4}
Suppose Assumptions \ref{a1}\mbox{-}\ref{a5} are satisfied, then as $n\to\infty$,
\begin{align}\label{(24)}
&\sigma_n(p_0, q^*)\notag\\
\leq&3\left(1+p_0c\right)\varepsilon n+\sum^{n-1}_{i=0}\left(w_{i+1}+\epsilon_{i+1}\right)^2+o\left(n\right), a.s.
\end{align}
\end{lemmax}

\emph{Proof}:
From (\ref{(4)}), (\ref{(5)}), (\ref{(7)}), (\ref{(13)}), (\ref{(14)}), (\ref{(22)}) we have
\begin{align}\label{(25)}
\sigma_n(p_0, q^*)=&\sum^{n-1}_{i=0}\left(\bar{\theta}^{\top}(p_0, q^*)\psi_i(p_0, q^*)+w_{i+1}+\epsilon_{i+1}\notag\right.\\
&\left.-\theta_n^{\top}(p_0, q^*)\psi_i(p_0, q^*)\right)^2.
\end{align}
So,
\begin{align}\label{(26)}
&\sigma_n(p_0, q^*)\notag\\
=&\tilde{\theta}^{\top}_n(p_0, q^*)\sum^{n-1}_{i=0}\psi_i(p_0, q^*)\psi^{\top}_i(p_0, q^*)\tilde{\theta}_n(p_0, q^*)\notag\\
&+2\tilde{\theta}^{\top}_n(p_0, q^*)\sum^{n-1}_{i=0}\psi_i(p_0, q^*)\left(w_{i+1}+\epsilon_{i+1}\right)\notag\\
&+\sum^{n-1}_{i=0}\left(w_{i+1}+\epsilon_{i+1}\right)^2.
\end{align}
From Theorem 1 of \citep{Jing2022} we get
\begin{align}\label{(27)}
&\tilde{\theta}^{\top}_n(p_0, q^*)\sum^{n-1}_{i=0}\psi_i(p_0, q^*)\psi^{\top}_i(p_0, q^*)\tilde{\theta}_n(p_0, q^*)\notag\\
\leq&\left(1+p_0c\right)\varepsilon n+o\left(n\right), a.s.,
\end{align}
and
\begin{align}\label{(28)}
&2\left|\tilde{\theta}^{\top}_n(p_0, q^*)\sum^{n-1}_{i=0}\psi_i(p_0, q^*)\left(w_{i+1}+\epsilon_{i+1}\right)\right|\notag\\
=&2\left|\tilde{\theta}^{\top}_n(p_0, q^*)\left(\bar{\theta}(p_0, q^*)-P^{-1}_n(p_0, q^*)\tilde{\theta}_n(p_0, q^*)\right)\right|\notag\\
\leq&2\left|\tilde{\theta}^{\top}_n(p_0, q^*)\bar{\theta}(p_0, q^*)\right|\notag\\
&+2\tilde{\theta}^{\top}_n(p_0, q^*)P^{-1}_n(p_0, q^*)\tilde{\theta}_n(p_0, q^*)\notag\\
\leq&2\left(1+p_0c\right)\varepsilon n+o\left(n\right), a.s.
\end{align}
From (\ref{(26)})\mbox{-}(\ref{(28)}) we obtain
\begin{align}\label{(29)}
&\sigma_n(p_0, q^*)\notag\\
\leq&3\left(1+p_0c\right)\varepsilon n+\sum^{n-1}_{i=0}\left(w_{i+1}+\epsilon_{i+1}\right)^2+o\left(n\right), a.s.
\end{align}
This completes the proof.\qed

Based on above lemmas, we can get the main theoretical result of the paper.
\begin{theorem}\label{t1}
Suppose Assumptions \ref{a1}\mbox{-}\ref{a5} are satisfied and $l_n$ satisfies
\begin{align}\label{(30)}
l_n\geq\left[5\left(1+p^*c\right)\varepsilon+\alpha_1\right]n,\quad\alpha_1>0,
\end{align}
and
\begin{align}\label{(31)}
&l_n\leq\frac{\alpha_2}{p^*}\left[a^2_{p_0}c_1-2\gamma\sqrt{c_3\left(1+p^*c\right)\varepsilon}-3\left(1+p^*c\right)\varepsilon\right]n,\notag\\
&0<\alpha_2<1,
\end{align}
then
\begin{align}\label{(32)}
\hat{p}_n\xrightarrow[n\rightarrow\infty]\ p_0, a.s.
\end{align}
\end{theorem}

\emph{Proof}:
First, we want to prove
\begin{align}\label{(33)}
\limsup_{n\to\infty}\hat{p}_n\leq p_0, a.s.
\end{align}

For $p>p_0$, similar with (\ref{(28)}) we have
\begin{align}\label{(34)}
&2\left|\tilde{\theta}^{\top}_n(p, q^*)\sum^{n-1}_{i=0}\psi_i(p, q^*)\left(w_{i+1}+\epsilon_{i+1}\right)\right|\notag\\
\leq&2\left(1+p^*c\right)\varepsilon n+o\left(n\right), a.s.
\end{align}
Similar with (\ref{(26)}) we have
\begin{align}\label{(35)}
&\sigma_n(p, q^*)\notag\\
=&\tilde{\theta}^{\top}_n(p, q^*)\sum^{n-1}_{i=0}\psi_i(p, q^*)\psi^{\top}_i(p, q^*)\tilde{\theta}_n(p, q^*)\notag\\
&+2\tilde{\theta}^{\top}_n(p, q^*)\sum^{n-1}_{i=0}\psi_i(p, q^*)\left(w_{i+1}+\epsilon_{i+1}\right)\notag\\
&+\sum^{n-1}_{i=0}\left(w_{i+1}+\epsilon_{i+1}\right)^2.
\end{align}
From (\ref{(34)}), (\ref{(35)}) we have
\begin{align}\label{(36)}
&\sigma_n(p, q^*)\notag\\
\geq&2\tilde{\theta}^{\top}_n(p, q^*)\sum^{n-1}_{i=0}\psi_i(p, q^*)\left(w_{i+1}+\epsilon_{i+1}\right)+\sum^{n-1}_{i=0}\left(w_{i+1}+\epsilon_{i+1}\right)^2\notag\\
\geq&-2\left(1+p^*c\right)\varepsilon n+\sum^{n-1}_{i=0}\left(w_{i+1}+\epsilon_{i+1}\right)^2+o\left(n\right), a.s.
\end{align}
From (\ref{(36)}) and Lemma \ref{l4} we have
\begin{align}\label{(37)}
&\sigma_n(p, q^*)-\sigma_n(p_0, q^*)\notag\\
\geq&-2\left(1+p^*c\right)\varepsilon n+\sum^{n-1}_{i=0}\left(w_{i+1}+\epsilon_{i+1}\right)^2+o\left(n\right)\notag\\
-&\left[3\left(1+p_0c\right)\varepsilon n+\sum^{n-1}_{i=0}\left(w_{i+1}+\epsilon_{i+1}\right)^2+o\left(n\right)\right]\notag\\
\geq&-5\left(1+p^*c\right)\varepsilon n+o\left(n\right), a.s.
\end{align}
From (\ref{(30)}) it can be seen that
\begin{align}\label{(38)}
l_n\cdot\left(p-p_0\right)\geq l_n\geq\left[5\left(1+p^*c\right)\varepsilon+\alpha_1\right].
\end{align}
From (\ref{(21)}), (\ref{(37)}), (\ref{(38)}) and noticing $\alpha_1>0$, we have
\begin{align}\label{(39)}
&\mathrm{min}_{\substack{p_0<p\leq p^*}}\left[L_n(p, q^*)-L_n(p_0, q^*)\right]\notag\\
\geq&-5\left(1+p^*c\right)\varepsilon n+l_n\cdot\left(p-p_0\right)+o\left(n\right)\notag\\
\geq&-5\left(1+p^*c\right)\varepsilon n+\left[5\left(1+p^*c\right)\varepsilon+\alpha_1\right]n+o\left(n\right)\notag\\
\geq&\alpha_1 n+o\left(n\right)\xrightarrow[n\rightarrow\infty]\ \infty, a.s.
\end{align}
So, (\ref{(33)}) is proved.

Next, we want to prove
\begin{align}\label{(40)}
\liminf_{n\to\infty}\hat{p}_n\geq p_0, a.s.
\end{align}
For $p<p_0$, from (\ref{(4)}), (\ref{(5)}), (\ref{(7)}), (\ref{(13)}), (\ref{(14)}), (\ref{(16)}), (\ref{(18)}), (\ref{(19)}) we have
\begin{align}\label{(41)}
&s_{i+1}-\theta^{\top}_n(p, q^*)\psi_i(p, q^*)\notag\\
=&s_{i+1}-\hat{\theta}^{\top}_n(p)\psi_i(p_0, q^*)\notag\\
=&\bar{\theta}^{\top}(p_0, q^*)\psi_i(p_0, q^*)+w_{i+1}+\epsilon_{i+1}-\hat{\theta}^{\top}_n(p)\psi_i(p_0, q^*)\notag\\
=&\tilde{\theta}^{\top}_n(p)\psi_i(p_0, q^*)+w_{i+1}+\epsilon_{i+1}.
\end{align}
From (\ref{(22)}), (\ref{(41)}) we have
\begin{align}\label{(42)}
&\sigma_n(p, q^*)\notag\\
=&\sum^{n-1}_{i=0}\left(s_{i+1}-\theta_n^{\top}(p, q^*)\psi_i(p, q^*)\right)^2\notag\\
=&\sum^{n-1}_{i=0}\left(\tilde{\theta}^{\top}_n(p)\psi_i(p_0, q^*)+w_{i+1}+\epsilon_{i+1}\right)^2\notag\\
=&\tilde{\theta}^{\top}_n(p)\sum^{n-1}_{i=0}\psi_i(p_0, q^*)\psi^{\top}_i(p_0, q^*)\tilde{\theta}_n(p)\notag\\
&+2\tilde{\theta}^{\top}_n(p)\sum^{n-1}_{i=0}\psi_i(p_0, q^*)\left(w_{i+1}+\epsilon_{i+1}\right)\notag\\
&+\sum^{n-1}_{i=0}\left(w_{i+1}+\epsilon_{i+1}\right)^2.
\end{align}
From (\ref{(13)}), (\ref{(18)}), (\ref{(19)}) we get
\begin{align}\label{(43)}
\left|\left|\tilde{\theta}^{\top}_n(p)\right|\right|^2\geq a^2_{p_0}>0.
\end{align}
From (\ref{(20)}), (\ref{(43)}) and Assumption \ref{a5} we have
\begin{align}\label{(44)}
&\tilde{\theta}^{\top}_n(p)\sum^{n-1}_{i=0}\psi_i(p_0, q^*)\psi^{\top}_i(p_0, q^*)\tilde{\theta}_n(p)\notag\\
=&\tilde{\theta}^{\top}_n(p)\left(\sum^{n-1}_{i=0}\psi_i(p_0, q^*)\psi^{\top}_i(p_0, q^*)+I-I\right)\tilde{\theta}_n(p)\notag\\
\geq&a^2_{p_0}\lambda^{(p_0, q^*)}_{min}(n-1)-\left|\left|\tilde{\theta}_n(p)\right|\right|^2\notag\\
\geq&a^2_{p_0}c_1n-\gamma^2.
\end{align}
\begin{align}\label{(45)}
&2\tilde{\theta}^{\top}_n(p)\sum^{n-1}_{i=0}\psi_i(p_0, q^*)\left(w_{i+1}+\epsilon_{i+1}\right)\notag\\
=&2\left|\left|\tilde{\theta}^{\top}_n(p)\left(\sum^{n-1}_{i=0}\psi_i(p_0, q^*)\psi^{\top}_i(p_0, q^*)+ I\right)^{\frac{1}{2}}\notag\right.\right.\\
&\left(\sum^{n-1}_{i=0}\psi_i(p_0, q^*)\psi^{\top}_i(p_0, q^*)+ I\right)^{-\frac{1}{2}}\notag\\
&\left.\left.\sum_{i=0}^{n-1}\psi_i(p_0, q^*)\left(w_{i+1}+\epsilon_{i+1}\right)\right|\right|.
\end{align}
From  Lemma \ref{l1}, Lemma \ref{l2}, Lemma \ref{l3} and (\ref{(45)}) we have
\begin{align}\label{(46)}
&2\tilde{\theta}^{\top}_n(p)\sum^{n-1}_{i=0}\psi_i(p_0, q^*)\left(w_{i+1}+\epsilon_{i+1}\right)\notag\\
\leq&2\left|\left|\tilde{\theta}^{\top}_n(p)\right|\right|\left|\left|\left(\sum^{n-1}_{i=0}\psi_i(p_0, q^*)\psi^{\top}_i(p_0, q^*)+ I\right)^{\frac{1}{2}}\right|\right|\notag\\
&\sqrt{\left(1+p_0c\right)\varepsilon n+o\left(n\right)}\notag\\
\leq&2\gamma\sqrt{\lambda^{(p_0, q^*)}_{max}(n-1)}\sqrt{\left(1+p_0c\right)\varepsilon n+o\left(n\right)}\notag\\
\leq&2\gamma\sqrt{c_3n}\sqrt{\left(1+p_0c\right)\varepsilon n+o\left(n\right)}\notag\\
=&2\gamma\sqrt{c_3\left(1+p_0c\right)\varepsilon+c_3o\left(1\right)}n\notag\\
\leq&2\gamma\sqrt{c_3\left(1+p_0c\right)\varepsilon}n+o\left(n\right).
\end{align}
From (\ref{(42)}), (\ref{(44)}) and (\ref{(46)}) it follows that
\begin{align}\label{(47)}
\sigma_n(p, q^*)\geq&a^2_{p_0}c_1n-2\gamma\sqrt{c_3\left(1+p_0c\right)\varepsilon}n\notag\\
&+\sum^{n-1}_{i=0}\left(w_{i+1}+\epsilon_{i+1}\right)^2+o\left(n\right).
\end{align}
From (\ref{(47)}) and Lemma \ref{a4} we have
\begin{align}\label{(48)}
&\sigma_n(p, q^*)-\sigma_n(p_0, q^*)\notag\\
\geq&a^2_{p_0}c_1n-2\gamma\sqrt{c_3\left(1+p_0c\right)\varepsilon}n+\sum^{n-1}_{i=0}\left(w_{i+1}+\epsilon_{i+1}\right)^2\notag\\
&-\left[3\left(1+p_0c\right)\varepsilon n+\sum^{n-1}_{i=0}\left(w_{i+1}+\epsilon_{i+1}\right)^2+o\left(n\right)\right]\notag\\
&+o\left(n\right)\notag\\
=&a^2_{p_0}c_1n-2\gamma\sqrt{c_3\left(1+p_0c\right)\varepsilon}n-3\left(1+p_0c\right)\varepsilon n+o\left(n\right)\notag\\
\geq&a^2_{p_0}c_1n-2\gamma\sqrt{c_3\left(1+p^*c\right)\varepsilon}n-3\left(1+p^*c\right)\varepsilon n\notag\\
&+o\left(n\right), a.s.
\end{align}
From (\ref{(31)}) it can be seen that
\begin{align}\label{(49)}
&l_n\cdot\left(p_0-p\right)\notag\\
\leq&\frac{\alpha_2}{p^*}\left[a^2_{p_0}c_1-2\gamma\sqrt{c_3\left(1+p^*c\right)\varepsilon}-3\left(1+p^*c\right)\varepsilon\right]n p^*\notag\\
=&\alpha_2\left[a^2_{p_0}c_1-2\gamma\sqrt{c_3\left(1+p^*c\right)\varepsilon}-3\left(1+p^*c\right)\varepsilon\right]n.
\end{align}
From (\ref{(21)}), (\ref{(48)}), (\ref{(49)}) and noticing $0<\alpha_2<1$, we have
\begin{align*}
&\mathrm{min}_{\substack{0\leq p<p_0}}\left[L_n(p, q^*)-L_n(p_0, q^*)\right]\notag\\
\geq&a^2_{p_0}c_1n-2\gamma\sqrt{c_3\left(1+p^*c\right)\varepsilon}n-3\left(1+p^*c\right)\varepsilon n+o\left(n\right)\notag\\
&-l_n\cdot\left(p_0-p\right)\notag\\
\geq&a^2_{p_0}c_1n-2\gamma\sqrt{c_3\left(1+p^*c\right)\varepsilon}n-3\left(1+p^*c\right)\varepsilon n+o\left(n\right)\notag\\
&-\alpha_2\left[a^2_{p_0}c_1-2\gamma\sqrt{c_3\left(1+p^*c\right)\varepsilon}-3\left(1+p^*c\right)\varepsilon\right]n\notag\\
=&\left(1-\alpha_2\right)\left[a^2_{p_0}c_1-2\gamma\sqrt{c_3\left(1+p^*c\right)\varepsilon}-3\left(1+p^*c\right)\varepsilon\right]n\notag\\
&+o\left(n\right)\notag\\
&\xrightarrow[n\rightarrow\infty]\ \infty, a.s.
\end{align*}
So, (\ref{(40)}) is proved.

From (\ref{(33)}), (\ref{(40)}) we know that
\begin{align}\label{(50)}
\hat{p}_n\xrightarrow[n\rightarrow\infty]\ p_0, a.s.
\end{align}
This completes the proof.\qed

\subsection{The estimation of $q_0$}
The quantized criterion $V_n(p, q)$ can be defined as
\begin{align}\label{(51)}
V_n(p, q):=\sigma_n(p, q)+v_n{\cdot}\left(p+q\right),
\end{align}
where
$
\sigma_n(p, q)
$
is defined in (\ref{(22)}) and the restrictions of $v_n$ will be given later.

The order estimation $\hat{q}_n$ of $q_0$ is defined as
\begin{align}\label{(52)}
\hat{q}_n:=\mathrm{argmin}_{\substack{0\leq q\leq q^*}}V_n(p^*, q).
\end{align}

\begin{lemmax}
Suppose Assumptions \ref{a1}\mbox{-}\ref{a2} are satisfied. Then, as $n\to\infty$, there is a constant $c_4>0$ such that
\begin{align}
\lambda^{(p^*, q_0)}_{max}(n)\leq c_4\left(n+1\right), a.s.,
\end{align}
where $\lambda^{(p^*, q_0)}_{max}(n)$ denotes the largest eigenvalue of
$
\sum_{i=0}^{n}\psi_i(p^*, q_0)\psi^{\top}_i(p^*, q_0)+I.
$
\end{lemmax}
\emph{Proof}: The proof is similar with Lemma \ref{l1}.\qed

\begin{lemmax}
Suppose Assumptions \ref{a1}\mbox{-}\ref{a4} and \ref{a6} are satisfied under the condition $q\leq q_0$, and define
\begin{align}
\hat{\theta}_n(q):=[&-a_{1n}(q), \cdots, -a_{p^*n}(q), b_{1n}(q), \cdots, b_{qn}(q),\notag\\
&\underbrace{0, \cdots, 0}_{q_0-q}]^{\top},
\end{align}
where $a_{in}(q)$, $b_{in}(q)$ are of $\theta_n(p^*, q)$.

Let
\begin{align}
\tilde{\theta}_n(q)=\bar{\theta}(p^*, q_0)-\hat{\theta}_n(q).
\end{align}
Then as $n\to\infty$, there is a constant $\gamma'$ such that
\begin{align}
\left|\left|\tilde{\theta}_n(q)\right|\right|\leq\gamma', a.s.
\end{align}
\end{lemmax}
\emph{Proof}: The proof is similar with Lemma \ref{l3}.\qed

\begin{theorem}\label{t2}
Suppose Assumptions \ref{a1}\mbox{-}\ref{a4} and \ref{a6} are satisfied and $v_n$ satisfies
\begin{align}\label{(53)}
v_n\geq\left[5\left(1+p^*c\right)\varepsilon+\beta_1\right]n,\quad\beta_1>0
\end{align}
and
\begin{align}\label{(54)}
&v_n\leq\frac{\beta_2}{q^*}\left[b^2_{q_0}c_2-2\gamma'\sqrt{c_4\left(1+p^*c\right)\varepsilon}-3\left(1+p^*c\right)\varepsilon\right]n,\notag\\
&0<\beta_2<1,
\end{align}
then
\begin{align}\label{(55)}
\hat{q}_n\xrightarrow[n\rightarrow\infty]\ q_0, a.s.
\end{align}
\end{theorem}

\emph{Proof}:
The proof is similar with Theorem \ref{t1}.\qed
\begin{remark}
By choosing suitable $\varepsilon$, $\alpha_1$, $\alpha_2$, $\beta_1$ and $\beta_2$ it can be made sure that
\begin{align*}
\Big[&5\left(1+p^*c\right)\varepsilon+\alpha_1,\notag\\
&\frac{\alpha_2}{p^*}\left(a^2_{p_0}c_1-2\gamma\sqrt{c_3\left(1+p^*c\right)\varepsilon}-3\left(1+p^*c\right)\varepsilon\right)\Big]
\end{align*}
and
\begin{align*}
\Big[&5\left(1+p^*c\right)\varepsilon+\beta_1,\notag\\
&\frac{\beta_2}{q^*}\left(b^2_{q_0}c_2-2\gamma'\sqrt{c_4\left(1+p^*c\right)\varepsilon}-3\left(1+p^*c\right)\varepsilon\right)\Big]
\end{align*}
are not empty sets. So, (\ref{(30)}), (\ref{(31)}), (\ref{(53)}) and (\ref{(54)}) are meaningful.
\end{remark}
\begin{remark}
\textcolor{blue}{Selecting $a_{p_0}$, $b_{q_0}$, $\gamma$ and $\gamma'$ in (\ref{(30)})\mbox{-}(\ref{(31)}), (\ref{(53)})\mbox{-}(\ref{(54)}) depends on the exact model and order of the system, and we do not have access to them. Actually, this limit is similar with the conditions in Theorem 7.1 of \citep{Chen1991}.}
\end{remark}
\section{Numerical example}
In this section, we will illustrate the theoretical result with a simulation example.

Consider the following ARX system:
$
y_{n}=a_1y_{n-1}+a_2y_{n-2}+b_1u_{n-1}+w_{n}, n=1, 2, ... ,
$
where the system noise $w_{n}$ follows $N(0, 1)$, $p_0=2$, $q_0=1$. $\theta=\left[a_1, a_2, b_1\right]^{\top}=\left[-0.7, -0.1, 1\right]^{\top}$. Let $y_n$ be quantized by (\ref{(3)}) under $\varepsilon=0.001$ and $\varepsilon=0.002$, $p^*=3$, $q^*=3$ and $p^*=6$, $q^*=6$, respectively.

With the selected \textcolor{blue}{$p\ (p\leq p^*)$} and \textcolor{blue}{$q\ (q\leq q^*)$}, we use the following algorithm to estimate $p_0$ and $q_0$.

\begin{algorithm}[h]
\renewcommand{\algorithmicrequire}{\textbf{Input:}}
\renewcommand{\algorithmicensure}{\textbf{Output:}}
\caption{\textcolor{blue}{The estimate of $p_0$ and $q_0$}}\label{Algorithm}
\textcolor{blue}{
\begin{algorithmic}[1]
\Require
$u_i$.
\Ensure
$\hat{p}_n$ and $\hat{q}_n$.
\State Compute $\theta_n(p, q)$ according to Eq. (\ref{(15)});
\State Compute $\sigma_n(p, q)$ according to Eq. (\ref{(22)});
\State Compute $L_n(p, q)$ according to Eq. (\ref{(21)});
\State Compute $V_n(p, q)$ according to Eq. (\ref{(51)});
\State Compute $\hat{p}_n$ according to Eq. (\ref{(23)});
\State Compute $\hat{q}_n$ according to Eq. (\ref{(52)}).
\end{algorithmic}}
\end{algorithm}

\textcolor{blue}{For the estimate of $p_0$, we chose $u_i$ to satisfy uniform distribution in $[-3, 3]$.} From (\ref{(30)})\mbox{-}(\ref{(31)}), when $\varepsilon=0.001$, let $l_n=0.006n$ and when $\varepsilon=0.002$, let $l_n=0.012n$. The trajectories of $\hat{p}_n$ are given by Fig. 1\mbox{-}4.

\textcolor{blue}{For the estimate of $q_0$, we chose $u_i$ to satisfy uniform distribution in $[-1, 1]$.} From (\ref{(53)})-(\ref{(54)}), when $\varepsilon=0.001$, let $v_n=0.006n$ and when $\varepsilon=0.002$, let $v_n=0.012n$. The trajectories of $\hat{q}_n$ are given by Fig. 5\mbox{-}8.

From Fig. 1\mbox{-}4, we can see that $\hat{p}_n$ converges to the true value $p_0$. From Fig. 5\mbox{-}8, we can see that $\hat{q}_n$ converges to the true value $q_0$. Moreover, the convergence rates of $\hat{p}_n$ and $\hat{q}_n$ are affected by the bounds \textcolor{blue}{$p^*$} and \textcolor{blue}{$q^*$}. To be more concrete, the larger the bounds, the slower convergence rates of $\hat{p}_n$ and $\hat{q}_n$.

\begin{figure}[!h]
\centering
\includegraphics[width=0.35\textwidth]{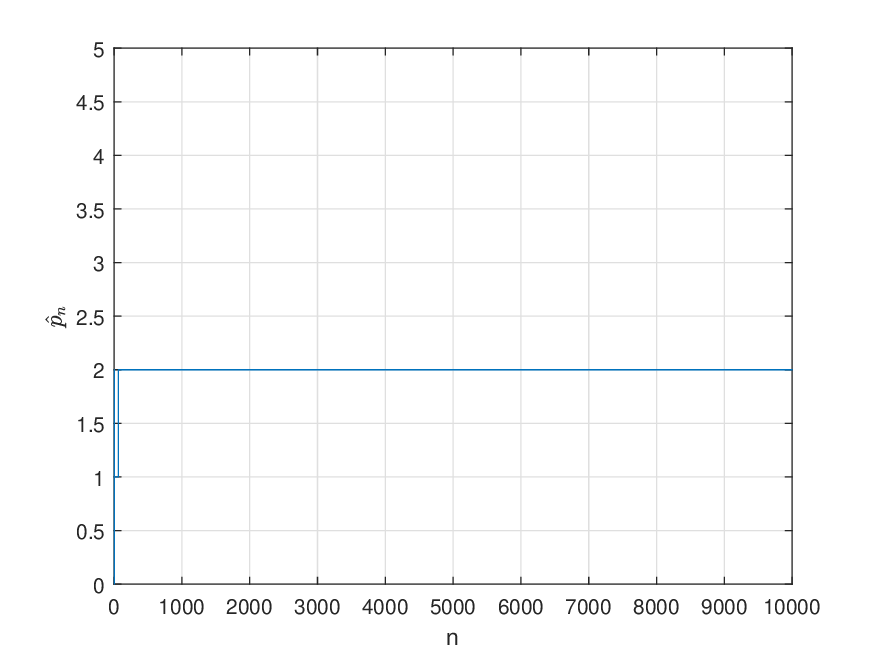}
\caption{The trajectories of $\hat{p}_n$ with $\varepsilon=0.001$, $p^*=3$, $q^*=3$}
\end{figure}

\begin{figure}[!h]
\centering
\includegraphics[width=0.35\textwidth]{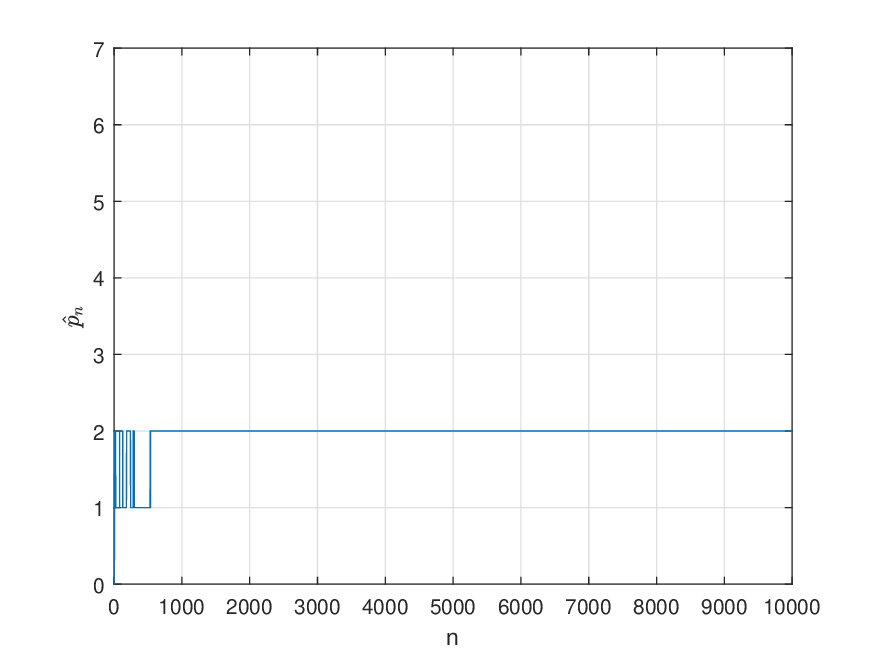}
\caption{The trajectories of $\hat{p}_n$ with $\varepsilon=0.001$, $p^*=6$, $q^*=6$}
\end{figure}

\begin{figure}[!h]
\centering
\includegraphics[width=0.35\textwidth]{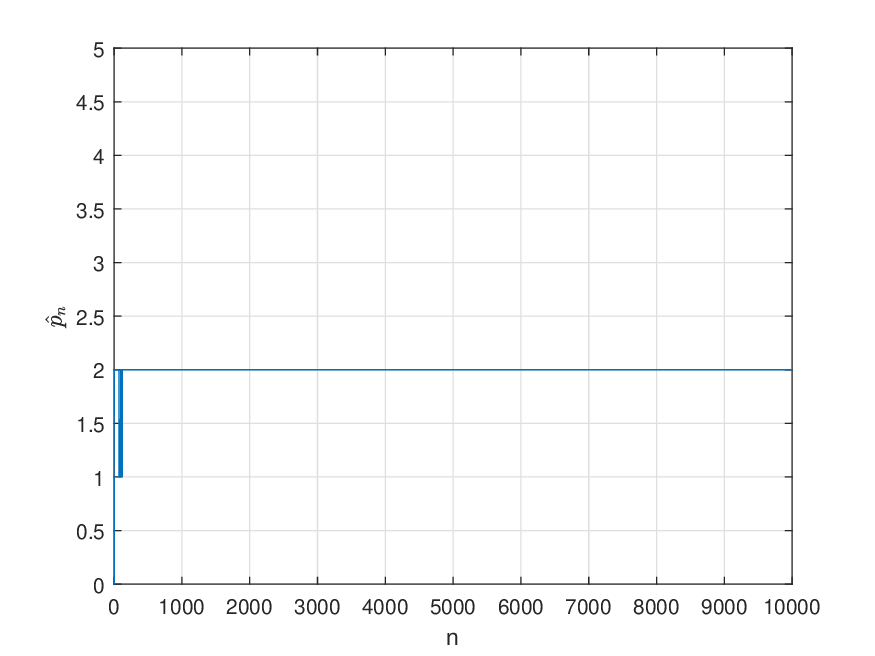}
\caption{The trajectories of $\hat{p}_n$ with $\varepsilon=0.002$, $p^*=3$, $q^*=3$}
\end{figure}

\begin{figure}[!h]
\centering
\includegraphics[width=0.35\textwidth]{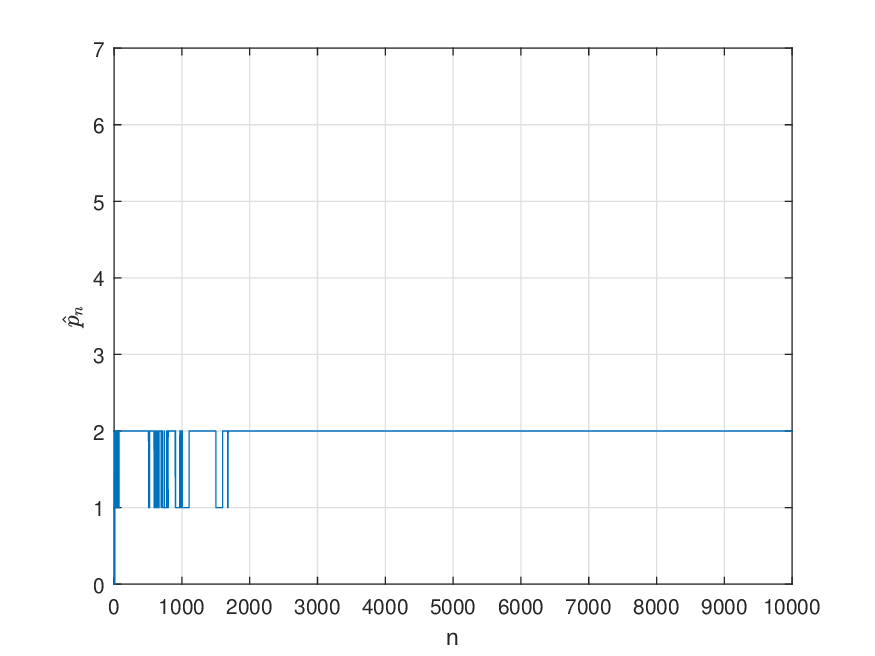}
\caption{The trajectories of $\hat{p}_n$ with $\varepsilon=0.002$, $p^*=6$, $q^*=6$}
\end{figure}

\begin{figure}[!h]
\centering
\includegraphics[width=0.35\textwidth]{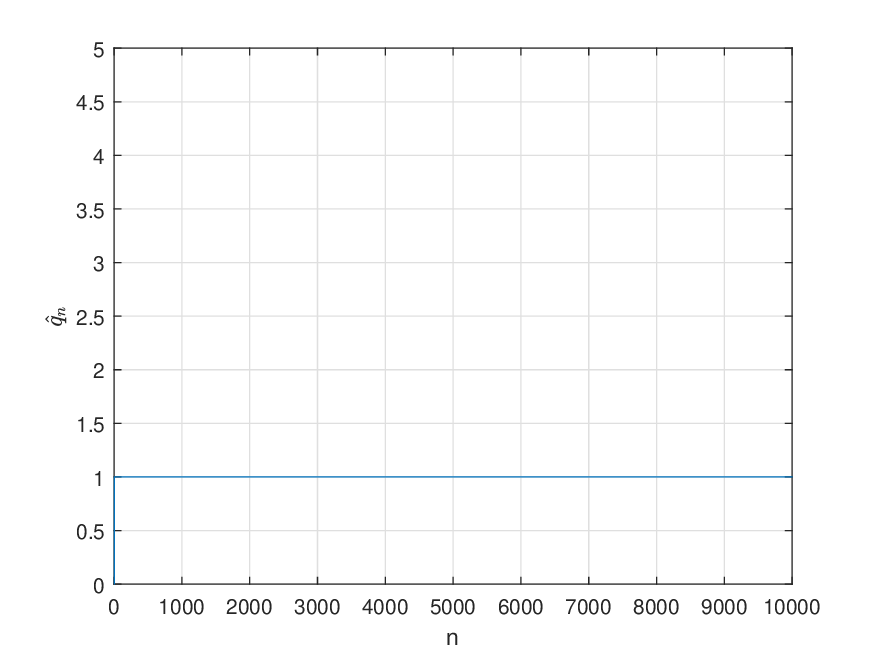}
\caption{The trajectories of $\hat{q}_n$ with $\varepsilon=0.001$, $p^*=3$, $q^*=3$}
\end{figure}

\begin{figure}[!h]
\centering
\includegraphics[width=0.35\textwidth]{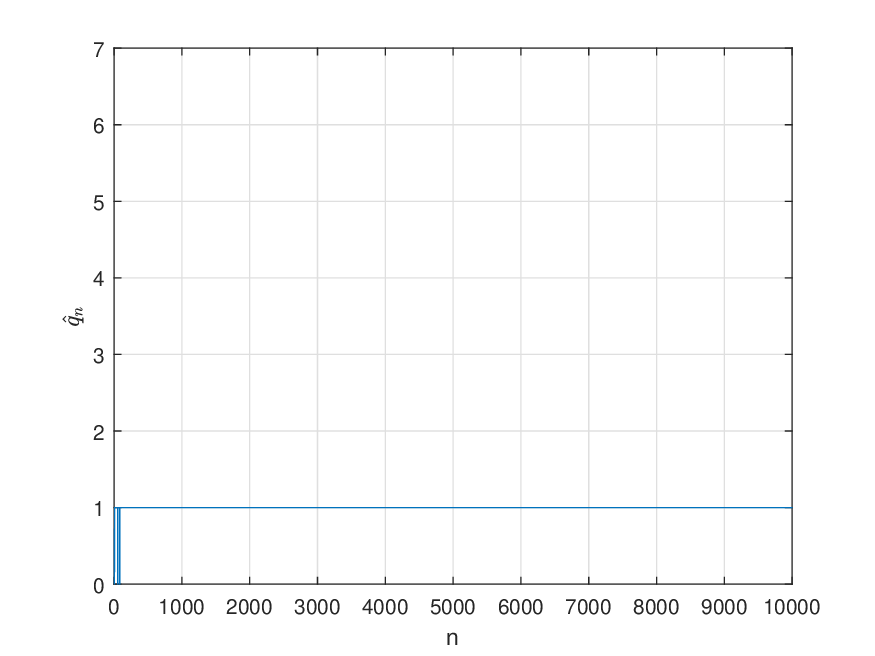}
\caption{The trajectories of $\hat{q}_n$ with $\varepsilon=0.001$, $p^*=6$, $q^*=6$}
\end{figure}

\begin{figure}[!h]
\centering
\includegraphics[width=0.35\textwidth]{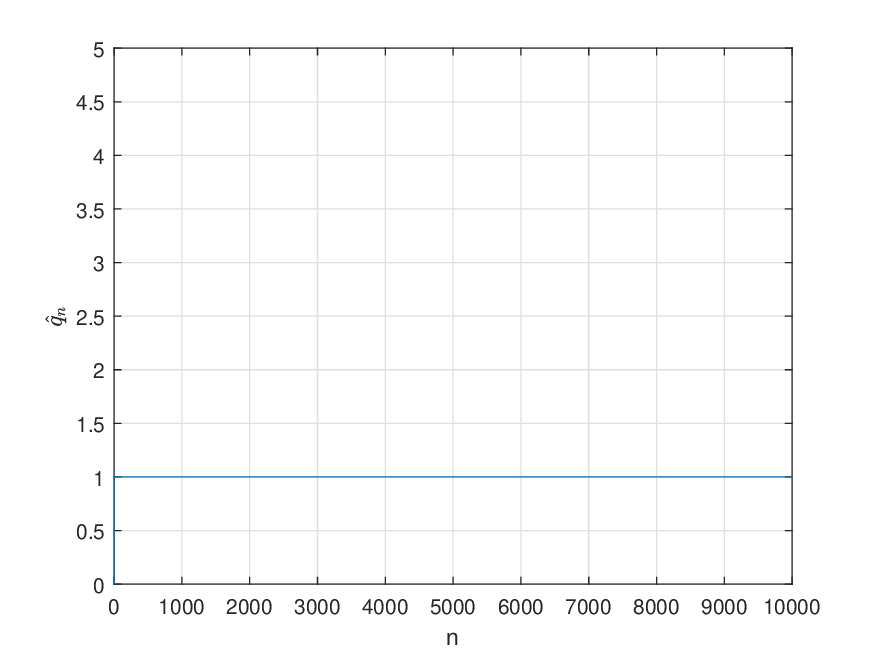}
\caption{The trajectories of $\hat{q}_n$ with $\varepsilon=0.002$, $p^*=3$, $q^*=3$}
\end{figure}

\begin{figure}[!h]
\centering
\includegraphics[width=0.35\textwidth]{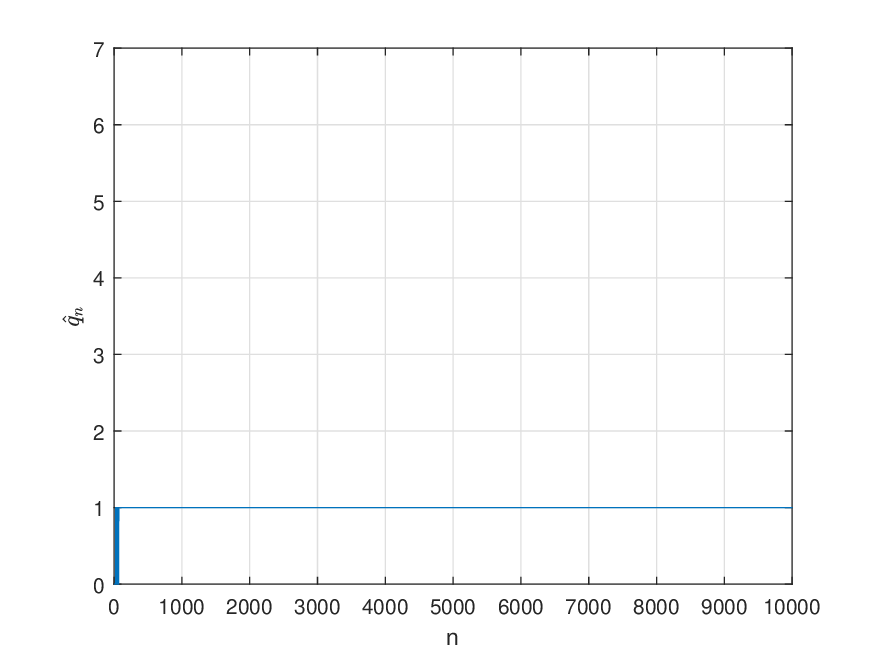}
\caption{The trajectories of $\hat{q}_n$ with $\varepsilon=0.002$, $p^*=6$, $q^*=6$}
\end{figure}

\section{Conclusion}
This paper considers the order estimation of ARX systems by using uniform quantized data. We design a novel criterion so as to estimate orders based on persistent excitation condition and some assumptions. Obviously, \citep{Jing2022} provides ideas for this paper and the least squares method is the key to the algorithm of this paper. It is shown that the estimated order is consistent. For further research, \textcolor{blue}{a method is
required for the verification of the assumptions and conditions introduced in Theorem \ref{t1} and \ref{t2}.} Another topic is how to reduce the amount of calculation. The methods proposed by \citep{Zhao2015} may be useful to solve such a problem.

\appendix
\section{Proof of Lemma \ref{l1}}\label{A}

From Assumptions \ref{a1}, \ref{a2}, Lemma B.3.3. of \citep{Goodwin1984}(Page 486) and law of large numbers, we know that there exists a positive constant $\hat{c}$ such that $\lim\limits_{n\to\infty}\frac{\sum_{i=0}^{n}y^2_i}{n+1}\leq\hat{c}$.

So, from (\ref{(3)}), (\ref{(7)}) and Assumption \ref{a1}, we know that
\begin{align*}
&\sum_{i=0}^{n}\left|\left|\psi_i(p_0, q^*)\right|\right|^2\\
=&\sum_{i=0}^{n}\left(s^2_i+s^2_{i-1}+\cdots+s^2_{i-p_0+1}\right.\\
&\left.+u^2_i+u^2_{i-1}+\cdots+u^2_{i-q^*+1}\right)\\
\leq&\sum_{i=0}^{n}2\left[y^2_i+\left(\frac{\varepsilon}{2}\right)^2+y^2_{i-1}+\left(\frac{\varepsilon}{2}\right)^2+\cdots\right.\\
&\left.+y^2_{i-p_0+1}+\left(\frac{\varepsilon}{2}\right)^2\right]+\sum_{i=0}^{n}\left(u^2_i+u^2_{i-1}+\cdots+u^2_{i-q^*+1}\right)\\
\leq&\left(2p_0\hat{c}+\frac{p_0\varepsilon^2}{2}+q^*\delta^2\right)\left(n+1\right).
\end{align*}

So, as $n\to\infty$, there exists a constant $c_3>0$ such that
\begin{align}\label{(A1)}
\lambda^{(p_0, q^*)}_{max}(n)=&\left|\left|\sum_{i=0}^{n}\psi_i(p_0, q^*)\psi^{\top}_i(p_0, q^*)+I\right|\right|\notag\\
\leq&\sum_{i=0}^{n}\left|\left|\psi_i(p_0, q^*)\psi^{\top}_i(p_0, q^*)\right|\right|+1\notag\\
=&\sum_{i=0}^{n}\left|\left|\psi_i(p_0, q^*)\right|\right|^2+1\notag\\
\leq&\left(2p_0\hat{c}+\frac{p_0\varepsilon^2}{2}+q^*\delta^2+1\right)\left(n+1\right), a.s.
\end{align}
Let $c_3=2p_0\hat{c}+\frac{p_0\varepsilon^2}{2}+q^*\delta^2+1$. This completes the proof.\qed

\section{Proof of Lemma \ref{l2}}\label{B}
From (\ref{(13)}) and (\ref{(14)}) we know that
\begin{align}\label{(B1)}
\bar{\theta}(p_0, q^*)=\left[-a_1, \cdots, -a_{p_0}, b_1, \cdots, b_{q_0}, 0, \cdots,0\right]^{\top},
\end{align}
and from (\ref{(4)}), (\ref{(5)}), (\ref{(7)}), (\ref{(13)}) and (\ref{(B1)}) it can be seen that
\begin{align}\label{(B2)}
s_{n+1}=\bar{\theta}^{\top}(p_0, q^*)\psi_n(p_0, q^*)+w_{n+1}+\epsilon_{n+1}.
\end{align}
From (\ref{(15)}) we know that
\begin{align}\label{(B3)}
\theta_n(p_0, q^*)=&\left(\sum^{n-1}_{i=0}\psi_i(p_0, q^*)\psi^{\top}_i(p_0, q^*)+I\right)^{-1}\notag\\
&\sum_{i=0}^{n-1}\psi_i(p_0, q^*)s_{i+1}\notag\\
=&P_n(p_0, q^*)\sum^{n-1}_{i=0}\psi_i(p_0, q^*)s_{i+1}.
\end{align}
From (\ref{(15)}) and (\ref{(B2)}) the estimated parameter error can be written as
\begin{align}\label{(B4)}
&\tilde{\theta}_n(p_0, q^*)\notag\\
=&\bar{\theta}(p_0, q^*)-\theta_n(p_0, q^*)\notag\\
=&\bar{\theta}(p_0, q^*)-P_n(p_0, q^*)\sum_{i=0}^{n-1}\psi_i(p_0, q^*)\notag\\
&\left(\psi_i^{\top}(p_0, q^*)\bar{\theta}(p_0, q^*)+w_{i+1}+\epsilon_{i+1}\right)\notag\\
=&\bar{\theta}(p_0, q^*)-P_n(p_0, q^*)\left(P^{-1}_n(p_0, q^*)-I\right)\bar{\theta}(p_0, q^*)\notag\\
&-P_n(p_0, q^*)\sum_{i=0}^{n-1}\psi_i(p_0, q^*)w_{i+1}\notag\\
&-P_n(p_0, q^*)\sum_{i=0}^{n-1}\psi_i(p_0, q^*)\epsilon_{i+1}\notag\\
=&P_n(p_0, q^*)\bar{\theta}(p_0, q^*)-P_n(p_0, q^*)\sum_{i=0}^{n-1}\psi_i(p_0, q^*)w_{i+1}\notag\\
&-P_n(p_0, q^*)\sum_{i=0}^{n-1}\psi_i(p_0, q^*)\epsilon_{i+1}.
\end{align}
From (\ref{(7)}) we have
\begin{align}\label{(B6)}
&\left|\left|\left(\sum^{n-1}_{i=0}\psi_i(p_0, q^*)\psi^{\top}_i(p_0, q^*)+ I\right)^{-\frac{1}{2}}\notag\right.\right.\\
&\left.\left.\sum_{i=0}^{n-1}\psi_i(p_0, q^*)\left(w_{i+1}+\epsilon_{i+1}\right)\right|\right|^2\notag\\
=&\left|\left|\left(\sum_{i=0}^{n-1}\psi_i(p_0, q^*)\left(w_{i+1}+\epsilon_{i+1}\right)\right)^{\top}P_n(p_0, q^*)\notag\right.\right.\\
&\left.\left.\left(\sum_{i=0}^{n-1}\psi_i(p_0, q^*)\left(w_{i+1}+\epsilon_{i+1}\right)\right)\right|\right|.
\end{align}
From Assumption \ref{a5} we know that
\begin{align}\label{(B7)}
\lambda^{(p_0, q^*)}_{min}(n)\geq c_1 n, a.s.
\end{align}

So, from (\ref{(B4)})\mbox{-}(\ref{(B7)}) and Theorem 1 of \citep{Jing2022} it can be seen that
\begin{align}\label{(B8)}
&\left|\left|\left(\sum^{n-1}_{i=0}\psi_i(p_0, q^*)\psi^{\top}_i(p_0, q^*)+I\right)^{-\frac{1}{2}}\notag\right.\right.\\
&\left.\left.\sum_{i=0}^{n-1}\psi_i(p_0, q^*)\left(w_{i+1}+\epsilon_{i+1}\right)\right|\right|^2\notag\\
=&\left|\left|\left(\bar{\theta}(p_0, q^*)-P^{-1}_n(p_0, q^*)\tilde{\theta}_n(p_0, q^*)\right)^{\top}P_n(p_0, q^*)\right.\right.\notag\\
&\left.\left.\left(\bar{\theta}(p_0, q^*)-P^{-1}_n(p_0, q^*)\tilde{\theta}_n(p_0, q^*)\right)\right|\right|\notag\\
\leq&2\tilde{\theta}^{\top}_n(p_0, q^*)P^{-1}_n(p_0, q^*)\tilde{\theta}_n(p_0, q^*)\notag\\
&+2\bar{\theta}^{\top}(p_0, q^*)P_n(p_0, q^*)\bar{\theta}(p_0, q^*)\notag\\
\leq&c'+\left(1+p_0c\right)\varepsilon n+O\left(\log n\right)+o\left(1\right)\notag\\
=&\left(1+p_0c\right)\varepsilon n+o\left(n\right), a.s.,
\end{align}
where $c'$ is a constant, and its definition can be found in Theorem 1 of \citep{Jing2022}.
This completes the proof.\qed

\begin{remark}
$2\tilde{\theta}^{\top}_n(p_0, q^*)P^{-1}_n(p_0, q^*)\tilde{\theta}_n(p_0, q^*)\leq c'+\left(1+p_0c\right)\varepsilon n+O\left(\log\lambda_{max}^{(p_0, q^*)}(n-1)\right)$ in (\ref{(B8)}) is similar with that in Theorem 1 (\citep{Jing2022}). To be more concrete, from (\ref{(13)}) and (\ref{(14)}), we just need to treat $0$ in (\ref{(13)}) as parameters to be estimated.
\end{remark}

\section{Proof of Lemma \ref{l3}}\label{C}
From (\ref{(15)}) we get
\begin{align}\label{(C2)}
&\left|\left|\theta_n(p, q^*)\right|\right|\notag\\
\leq&\left|\left|\left(\sum^{n-1}_{i=0}\psi_i(p, q^*)\psi^{\top}_i(p, q^*)+I\right)^{-1}\right|\right|\notag\\
&\left|\left|\sum^{n-1}_{i=0}\psi_i(p, q^*)s_{i+1}\right|\right|.
\end{align}
From Assumption \ref{a5} it can be seen that
\begin{align}\label{(C3)}
\left|\left|\left(\sum^{n-1}_{i=0}\psi_i(p, q^*)\psi^{\top}_i(p, q^*)+I\right)^{-1}\right|\right|\leq\frac{1}{c_1n}, a.s.,
\end{align}
By (\ref{(C2)}), (\ref{(C3)}) and Assumptions \ref{a1}, \ref{a2} it can be seen that
$
\left|\left|\theta_n(p, q^*)\right|\right|
$
is bounded (a.s.).

From (\ref{(13)}), (\ref{(14)}) and Assumption \ref{a3} we know that $\left|\left|\bar{\theta}(p_0, q^*)\right|\right|$ is bounded.

So, there is a constant $\gamma$ such that
\begin{align}\label{(C4)}
\left|\left|\tilde{\theta}_n(p)\right|\right|=&\left|\left|\bar{\theta}(p_0, q^*)-\hat{\theta}_n(p)\right|\right|\notag\\
\leq&\left|\left|\bar{\theta}(p_0, q^*)\right|\right|+\left|\left|\theta_n(p, q^*)\right|\right|\notag\\
\leq&\gamma, a.s.
\end{align}
This completes the proof.\qed

\end{document}